\documentclass{article}
\usepackage{graphicx}
\usepackage{todonotes}

\begin{document}

\font\bbb=msbm10
\def\R{\hbox{\bbb\char82}}
\def\Z{\hbox{\bbb\char90}}
\def\nbd{\hbox{\scriptsize{nbd}}}
\newcommand{\note}[1]{\todo[inline, backgroundcolor=lightgray]{ NOTE: #1}}

\newcommand{\J}{{\cal J}}

\setlength{\parindent}{0cm}
\setlength{\parskip}{3mm plus1mm minus1mm}

\title{Revisiting Comins, Hassell, and May}

\author{Stewart D. Johnson}

\date{\today \\ \tiny \copyright\ Stewart Johnson} 

\maketitle

\begin{abstract}

\noindent
In a sequence of papers in the early 1990's, Comins, Hassell, and May investigated Nicholson-Baily dynamics in a spatial implementation with diffusion. They delineated four types of dynamic behavior:  {\sl crystalline lattices, spatial chaos, spirals,}  and {\sl hard-to-start spirals.} We revisit their results with current computational methods, and develop a sampling technique to estimate the Lyapunov spectrum. We find that  {\sl spatial chaos, spirals,}  and {\sl hard-to-start spirals} are not separate categories, but part of a spectrum of behavior. We more thoroughly investigate the crystalline structures. We show that the Lyapunov sampling method can be used to find bifurcation curves in parameter space, and demonstrate an interesting spatial chimera. 
\smallskip


\noindent
Mathematics Subject Classification, primary: 37-04, 37M25; secondary: 37G, 34C15, 34C23, 34C28, 34D08, 58J55 

\end{abstract}

\section{Spatial Nicholson-Bailey}\label{sec:SNB}

\subsection{The Model}

The Nicholson-Bailey dynamic was introduced to model the interactions between a host/prey population of size $x$ and a parasite/predator population of size $y$, both taken as non-negative real values \cite{Nich33, Nich35}. Using discrete time $t\in\Z$, the dynamic is defined by:
\renewcommand{\arraystretch}{1.4 }
\begin{equation}\label{eqn:nb}
\begin{array}{rcccl}
x^{t+1} &=& f(x^t,y^t) &=& \lambda\,x^t\,e^{-y^t}\\
y^{t+1} &=& g(x^t,y^t) &=& x^t\,(1-e^{-y^t})
\end{array}
\end{equation}
The parameter $\lambda$ is the normalized growth rate, and both populations rapidly collapse for $0\le\lambda\le 1$. For $\lambda>1$ the system has a single coexistent fixed point at
\begin{equation}\label{eqn:fp}
\begin{array}{rcl}
x^* &=& \frac{\lambda\ln\lambda}{\lambda-1}\\
y^* &=& \ln\lambda
\end{array}
\end{equation}
which is repelling with complex eigenvalues. All positive initial values other than the fixed point lead to oscillations that increase without bound. The original biological application is for systems that undergo increasing oscillations leading to complete collapse (\cite{Nich33, Nich35}).

Though the individual system is highly unstable, a spatial implementation with diffusion can produce bounded population dynamics with a wide variety of patterns and behaviors.  Specifically, consider a rectangular integer lattice of population cells $(x_{m,n},y_{m,n})$ for $0 \le n\le N$ and $0\le m \le M$. The Nicholson-Bailey dynamic is applied to each cell
$$
\begin{array}{rcl}
\tilde{x}_{m,n}^{t}  &=& f(x_{m,n}^t,y_{m,n}^t)\\
\tilde{y}_{m,n}^{t}  &=& g(x_{m,n}^t,y_{m,n}^t)
\end{array}
$$
and followed by diffusion with migration rates $\mu_x$ and $\mu_y$:
$$
\begin{array}{rcl}
x_{m,n}^{t+1}  &=& (1-\mu_x)\,\tilde{x}_{m,n}^t + \frac{\mu_x}{|\nbd(m,n)|}\,\sum_{(i,j)\in\nbd(m,n)} \tilde{x}_{i,j}^t \\
y_{m,n}^{t+1}  &=& (1-\mu_y)\,\tilde{y}_{m,n}^t + \frac{\mu_y}{|\nbd(m,n)|}\,\sum_{(i,j)\in\nbd(m,n)} \tilde{y}_{j,i}^t
\end{array}
$$
Here $\hbox{nbd}(m,n)$ denotes the set of cells neighboring the cell at $(m,n)$. Typical choices are the standard $4$ and $8$ cell neighborhoods, and sometimes the  $6$ cell neighborhood using a hexagonal grid. The dynamics are remarkable similar across these three choices, and we standardize on an 8-cell neighborhood. 

Typical boundary conditions include absorbing, reflecting, or toroidal (wrap-around). For $M,N$ sufficiently large (around $90$ or so) the boundary assumptions typically  have little effect on the overall dynamics for homogeneous cases. 

\subsection{Comins, Hassell, and May}\label{sec:CHM}

Comins, Hassell, and May presented the above model in their 1991 paper \cite{CHM91}, with supplementary research and descriptions appearing in \cite{CHM92} and \cite{May95}.

They employed a $30\times 30$ grid with reflective boundaries, considered two growth rates $\lambda=2$ and $\lambda=10$, and explored the dynamics for varying migration rates $\mu_x$ and $\mu_y$.

\begin{figure}[ht]
\begin{center}
  \includegraphics[width=.7\linewidth]{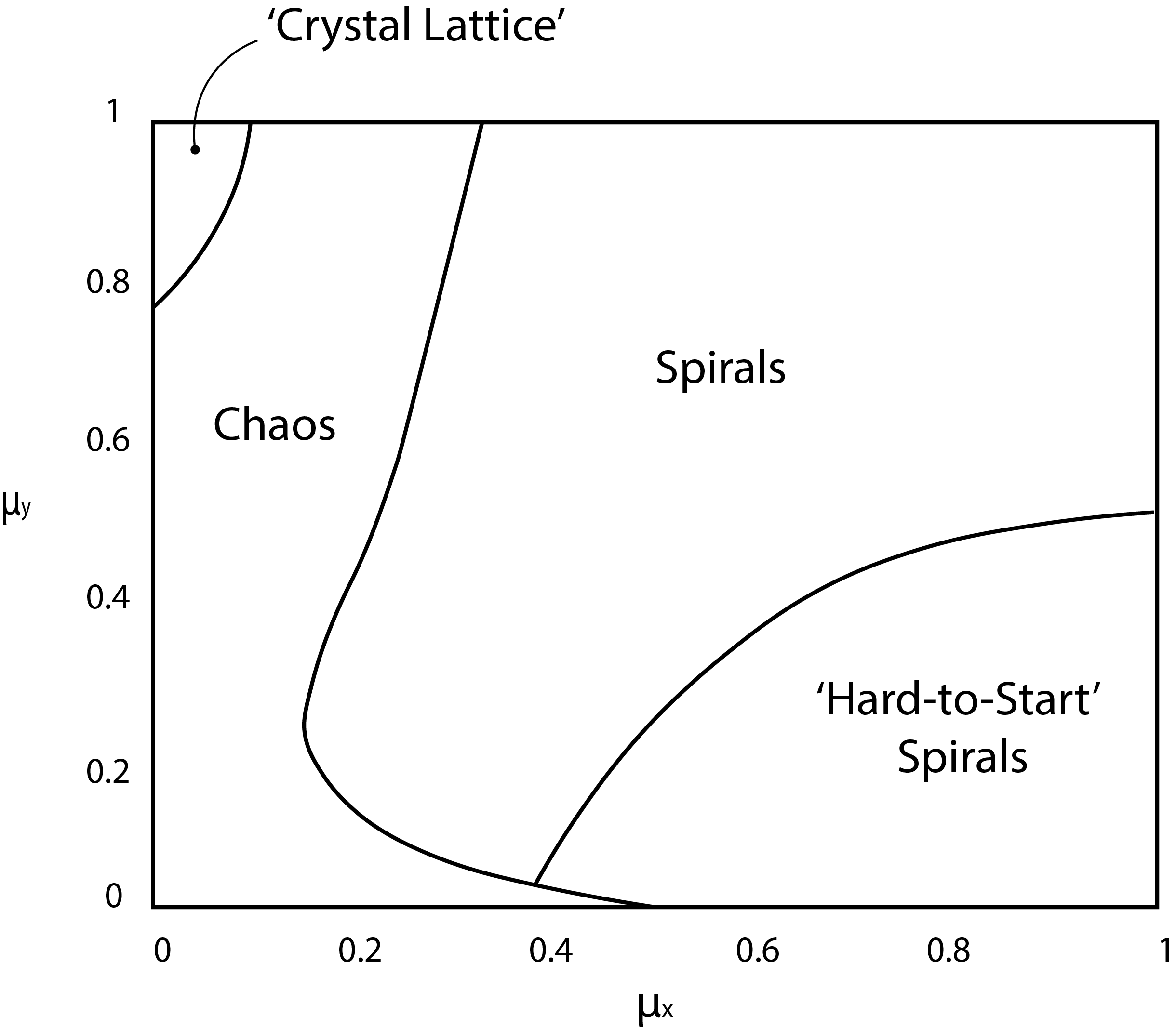}
\caption{Comins, Hassell, and May Diagram.}
\label{fig:CHM}
\end{center}
\end{figure}

They noted three distinct types of emergent patterns: {\sl crystal lattice, spatial chaos, spirals,}  and {\sl hard-to-start spirals,} and  mapped out regions in the $\mu_x \times \mu_y$ space that correspond to each type of pattern (see figure \ref{fig:CHM}). The boundaries of these regions were imprecise and subjective, and the assessment of chaos was by visual inspection. The following diagram appears in a number of publications \cite{CHM92,CHM94,SpatEco,OttarJordi}.

\subsection{Revisiting}

Revisiting Comins, Hassell, and May with modern methods, we find the following main conclusions:

\begin{itemize}
\item{ A serious limitation of the Comins, Hassell, and May research was the grid size of $30\times30$. The patterns and emergent structures inherent in these systems are on approximately this scale, and the dynamics are thus confounded by boundary conditions. Specifically, the {\sl hard-to-start spirals} case disappears in larger grids as the spirals have more room to form and stabilize.}

\item{There is no sharp transition between spirals and chaos, but more of a gradual change in the presence of spirals, waves, and more chaotic behavior. Most of the parameter region is chaotic, with smooth variation in the proportion of positive Lyapunov exponents and little variation in the maximum exponent. This is explored in section \ref{sec:chaos}}

\item{Notable behaviors appear at the fringes of the parameter region, including the {\sl crystal lattice} corner for small $\mu_x$ and large $\mu_y$ identified by Comins, Hassell, and May. Here we find that all Lyapunov exponents become negative as a stable global crystal lattice takes hold. This is explored in section \ref{sec:crystal}}

\end{itemize}

\section{Lyapunov Exponents}\label{sec:lyap}

\subsection{Computing Lyapunov Exponents}\label{sec:CompLCE}

The notion of chaos persists in being enigmatic. For low dimensional systems, estimating the Lyapunov Characteristic Exponents (LCE's) is a well established method for determining the presence or absence of chaotic behavior. In 1994, David Rand proposed that for high-dimensional spatial systems, chaos would be characterized by a positive lower bound for the proportion of positive LCE's as the dimension of the system increases \cite{Rand94}.  He demonstrated with a computation of LCE's for a one-dimensional coupled map lattice of  Nicholson-Bailey cells of length $128$.

Numerical analysis of LCE's is computationally expensive. Some earlier attempts at using Lyapunov computations to explore chaos in two-dimensional lattices reduced to time-series analysis of overall population size \cite{SoleValls,WakanoHauert}.

Full analysis of the Lyapunov spectrum involves computing the product of Jacobians of the entire system over a sequence of iterates. To avoid instabilities in the product, each Jacobian needs to be orthogonalized. Methodology is specified by Hubertus von Bremen in \cite{Hub1}, and a brief sketch of the theory is given in the appendix I.

Grids for exploring spatial systems need to be significantly larger that emergent structures to minimize boundary effects, and sufficiently large to reasonably address Rand's criteria for chaos which requires a lower bound on the proportion of positive LCE's as the size of the system expands. In the Nicholson-Bailey system, the state of each cell is given by two variables, and so the full Jacobian for an $N\times M$ grid is a $2MN \times 2MN $ matrix. For a $256\times256$ system this would involve a Jacobians of size $131,072\times 131,072$. Thousands of iterates are needed for reasonable convergence to the Lyapunov spectrum. Orthogonal computations on this scale are impractical.

\subsection{Subgrid Sampling}\label{sec:Samp}

For many planar systems the dynamics are spatially homogeneous. It is therefore possible to approximate LCE's using subgrids. Specifically, for a large $N\times M$ system, one computes the Lyapunov spectrum on a sampling of smaller $n\times m$ ($n<N$, $m<M$) subgrids. For each subgrid, this is equivalent to computing the Lyapunov exponents for the $n\times m$ system where the boundary conditions are determined by the larger system.   

The maximum LCE of a subgrid is necessarily smaller than that of the larger system, but proves to be a good approximate for a sufficiently large subgrid. Average LCE and proportion of positive LCE's are similarly well approximated. A numerical demonstration of the efficacy of subgrid sampling is given is appendix II.

\section{Chaotic Region}\label{sec:chaos}

To explore parameter space of spatial Nicholson-Bailey using subgrid sampling we standardize on $\lambda=2$ and compute the LCE's over a $20\times 20$ grid of parameter values for $0< \mu_x < 1$ and $0< \mu_y < 1$. For each $(\mu_x,\mu_y)$ pair, we consider a $768\times512$ lattice with wrapped boundaries, and seed each cell with uniformly distributed perturbation ($\pm 0.05$) of the fixed point $(x^*,y^*)$ from (\ref{eqn:fp}). The system is relaxed with $1,000,000$ initial iterates. We select three disjoint $32\times 32$ subgrids, and for the next $6,000$ iterates, the Jacobian of each subgrid is computed (a $2048 \times 2048$ matrix), QR-factored, and the LCE's are accumulated following methods as described in \cite{Hub1} and Appendix I.

\begin{figure}[!htbp]
\begin{center}
  \includegraphics[width=.65\linewidth]{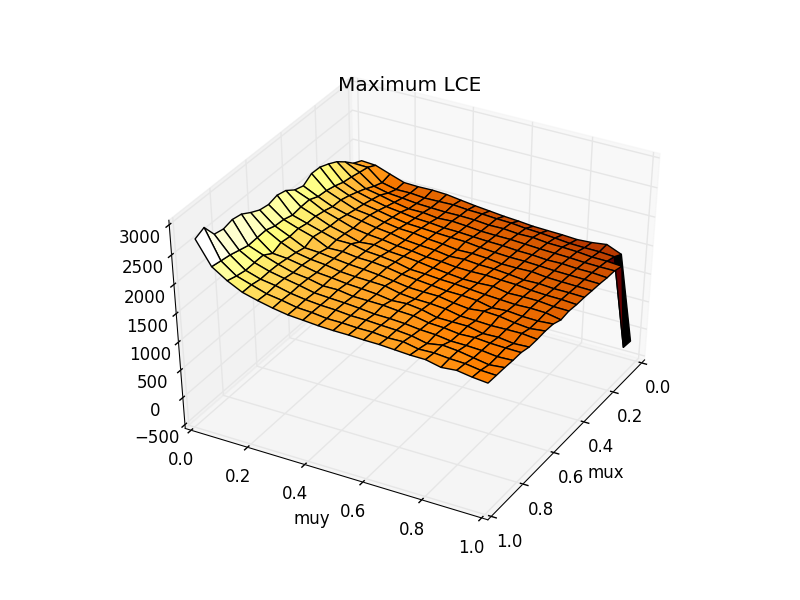}
  %\hskip6pt
  \includegraphics[width=.65\linewidth]{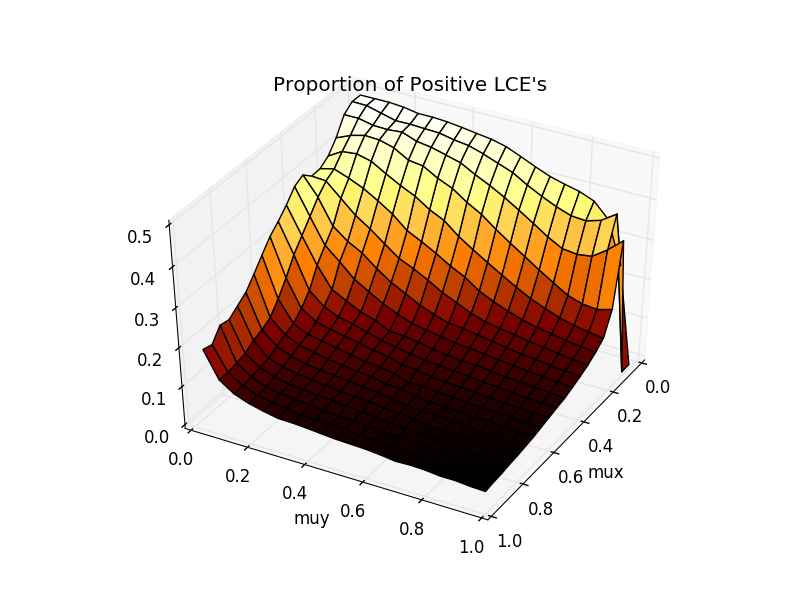}
  %\hskip6pt
  \includegraphics[width=.65\linewidth]{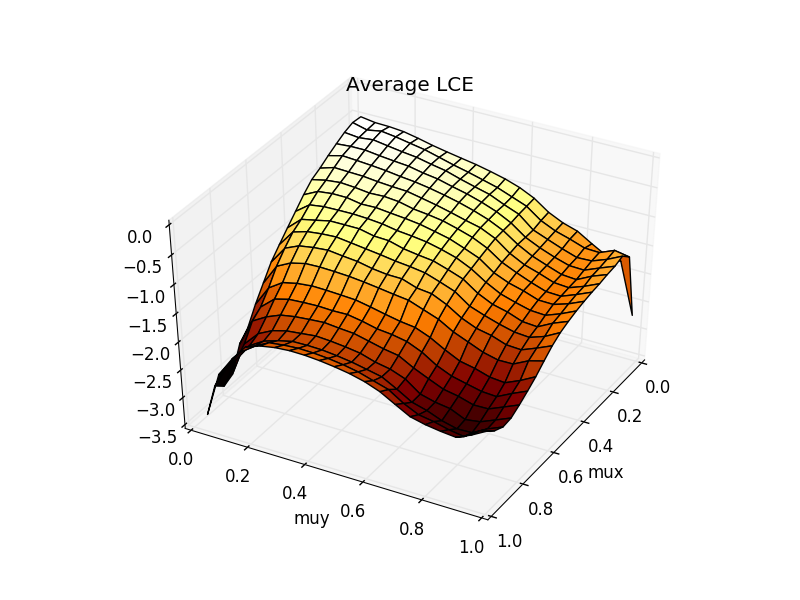}
\caption{Values of maximum LCE, proportion of positive LCE's, and average LCE for $\lambda=2$ and a grid of $\mu_x$, $\mu_y$ values.}
\label{fig:surf}
\end{center}
\end{figure}

Values of the maximum LCE show a remarkable consistency over $\mu_x$ and $\mu_y$ except for small values of $\mu_y$, and the \sl crystal lattice \rm corner of small $\mu_x$ and large $\mu_y$. The proportion of positive LCE's shows more variation, with a dramatic drop for  small $\mu_x$ and large $\mu_y$. Average LCE shows significant drops for both small $\mu_x$, large $\mu_y$ and small $\mu_y$, large $\mu_x$ (see figure \ref{fig:surf}).

\begin{figure}[!htbp]
\begin{center}
\includegraphics[width=.6\linewidth]{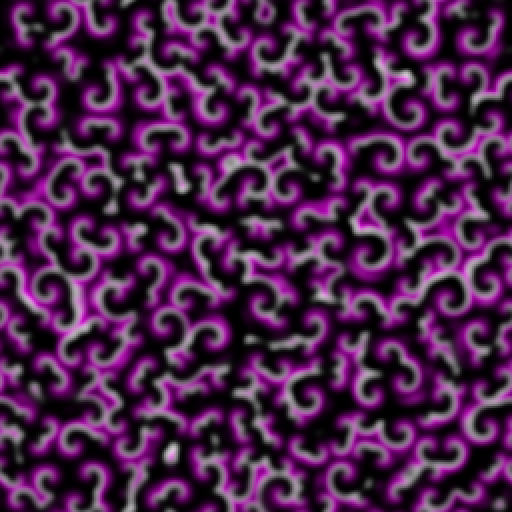}
\vskip10pt
\includegraphics[width=.6\linewidth]{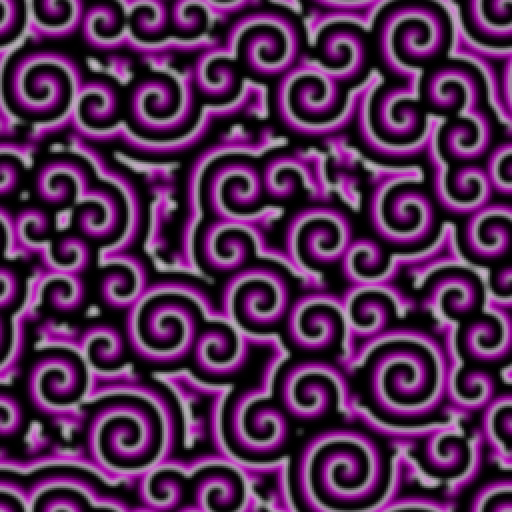}
\caption{Examples of large scale behavior for $\lambda=2$ on a $256 \times 256$ grid; (a) chaotic behavior for $\mu_x=0.6$, $\mu_y=0.8$, (b) BZ spirals for $\mu_x=0.9$, $\mu_y=0.5$. Purple indicates high $x$ values and gray indicates high $y$ values. }
\label{fig:examp}
\end{center}
\end{figure}

While there are certainly parameter values that lead to characteristic BZ spirals and others that seem highly chaotic with no discernible spiral patterns (see figure \ref{fig:examp}), there is no sharp transition between these behaviors and nothing in the pattern of LCE's would distinguish between cases of {\sl spirals} and {\sl chaos} cited by Comins, Hassell, and May.

\section{Crystal Lattice Corner}

Of interest is the corner of small $\mu_x$ and large $\mu_y$, corresponding to low levels of host/prey diffusion and high levels of parasite/predator diffusion, which gives rise to crystalline patterns identified by Comins, Hassell, and May (figure \ref{fig:CHM}). 

\subsection{Lattices, Islands, and Waves }\label{sec:crystal}

\begin{figure}[!htbp]
\begin{center}
  \includegraphics[width=.55\linewidth]{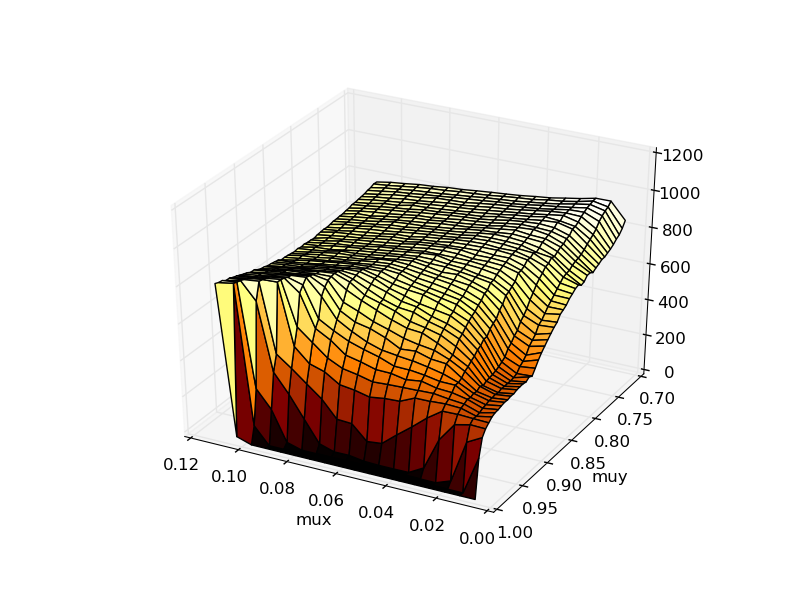}
  %\hskip6pt
  \includegraphics[width=.4\linewidth, height=.6\linewidth]{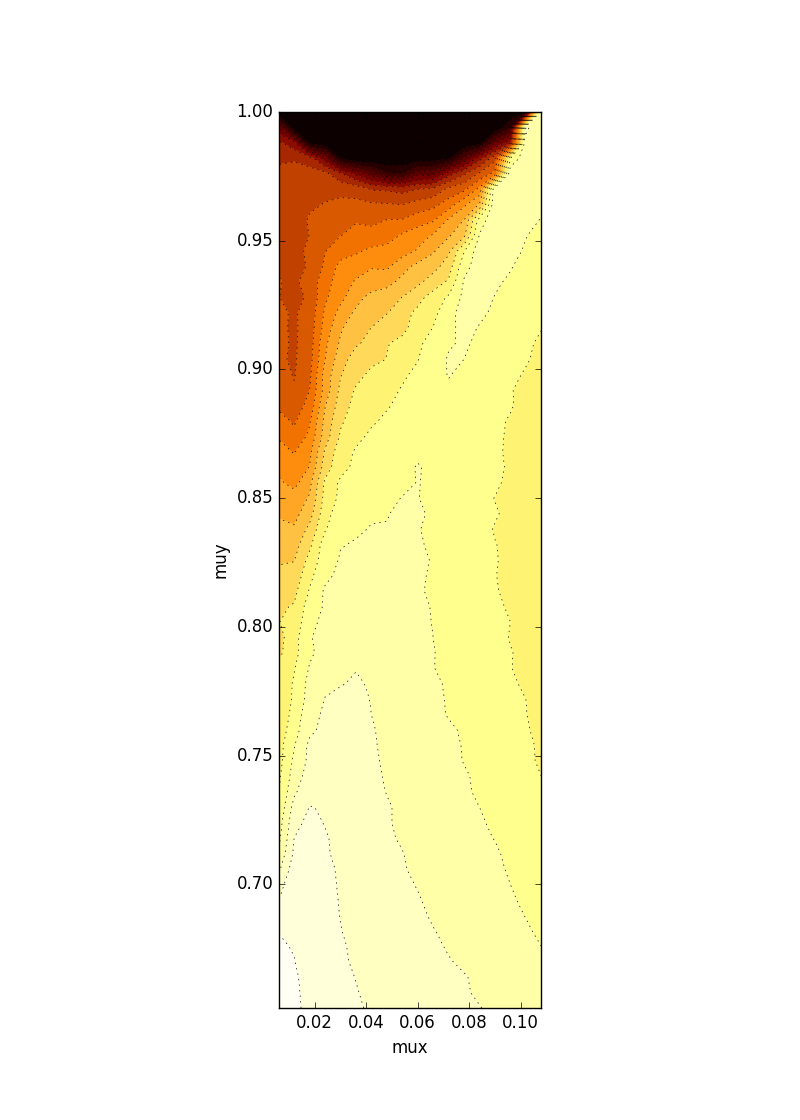}
\caption{Proportion of positive LCE's for $\lambda=2$, small $\mu_x$, and large $\mu_y$.}
\label{fig:NW}
\end{center}
\end{figure}

Investigating this corner with subgrid sampling reveals a semi-circular region near the $\mu_y=1$ boundary (the dark region in figure \ref{fig:NW}) which has no positive Lyapunov exponents, and in which the system settles into a stable global crystal lattice.  

\begin{figure}[!htbp]\label{fig:crystals}
\begin{center}
  \includegraphics[width=.24\linewidth]{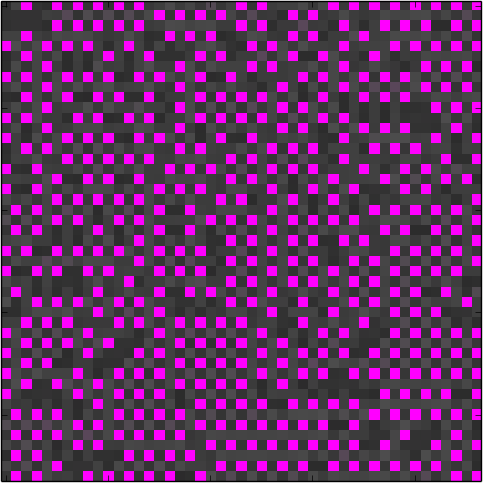}
  \includegraphics[width=.24\linewidth]{crysA.png}
  \includegraphics[width=.24\linewidth]{crysA.png}
  \includegraphics[width=.24\linewidth]{crysA.png}
  \vskip 10pt
  \includegraphics[width=.24\linewidth]{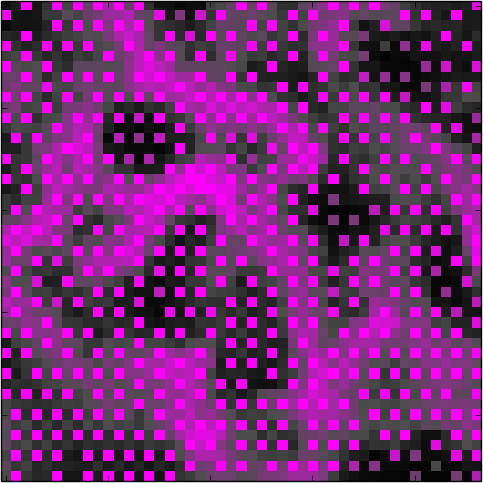}
  \includegraphics[width=.24\linewidth]{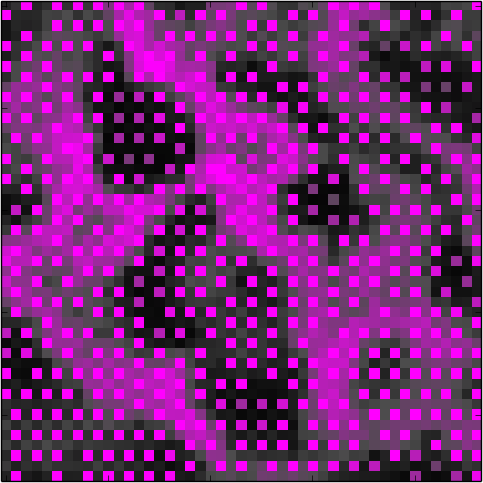}
  \includegraphics[width=.24\linewidth]{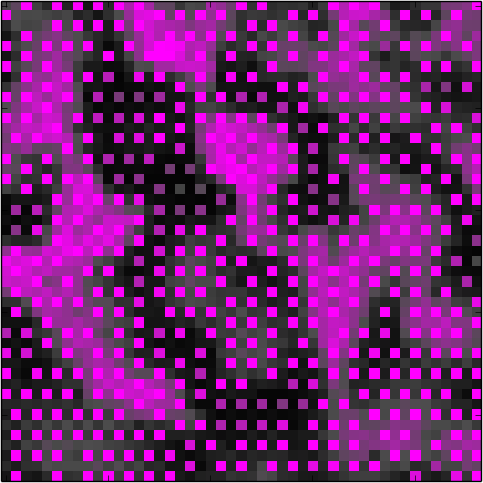}
  \includegraphics[width=.24\linewidth]{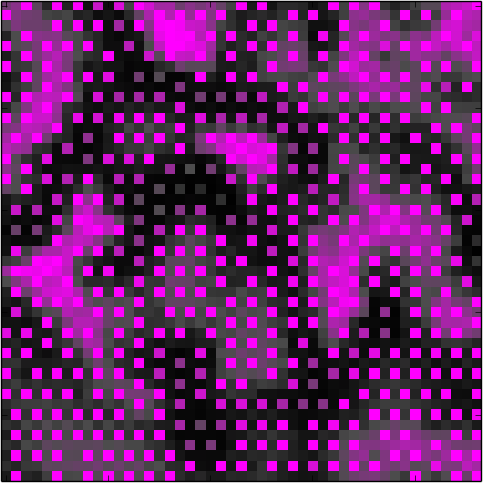}
  \vskip 10pt
  \includegraphics[width=.24\linewidth]{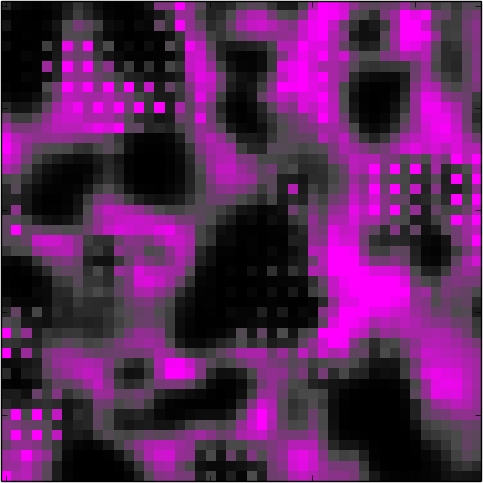}
  \includegraphics[width=.24\linewidth]{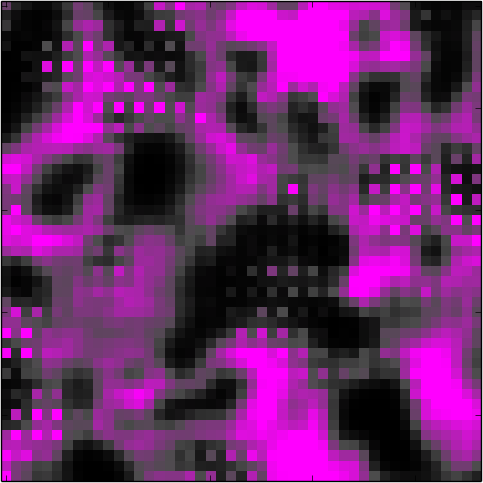}
  \includegraphics[width=.24\linewidth]{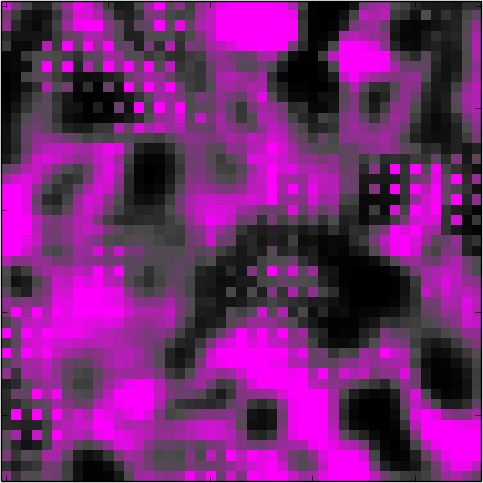}
  \includegraphics[width=.24\linewidth]{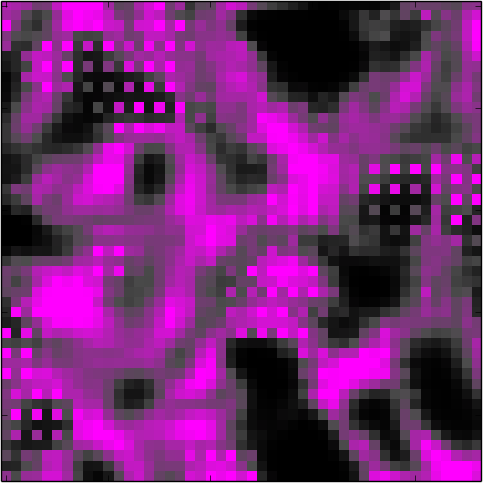}
  \vskip 10pt
  
\caption{Crystalline structures for $\lambda=2$ and $\mu_x=.08$, over $3$ iterations on a $32\times 32$ grid. Top row: $\mu_y=0.99$ fixed crystal lattice, no positive LCE's. Middle row: $\mu_y=0.98$ stable crystal lattice with a background of chaotic waves. Bottom row: $\mu_y=0.95$, unstable disconnected pockets of crystal structure in a sea of chaotic waves.  }
\end{center}
\end{figure}

A larger triangular region in the \sl crystal lattice \rm corner produces crystalline patterns, manifest as either stable global patterns with background waves, or transient islands of crystalline patterns that are destabilized by chaotic waves. Samples are shown in figure \ref{fig:crystals}. This region is related to a sharp decline in the proportion of positive LCE's; closely following the saddle ridge shown in the contour plot in figure \ref{fig:NW}, and duplicating the region found by  Comins, Hassell, and May (figure \ref{fig:CHM}). 

The crystal lattice is formed by cells with high $x$ values and low $y$ values, surrounded by eight neighboring cells with low $x$ values and high $y$ values. For $\mu_y \approx 1$ and $0<\mu_x<0.10$ these elements arrange themselves in a global stable pattern demonstrated in the top of figure \ref{fig:crystals}. The term {\sl crystal lattice} is well chosen, as large scale patterns will exhibit grain boundaries where the crystal pattern is mismatched. 

For a semi-circular region near the $\mu_y=1$ boundary (dark area in figure \ref{fig:NW}), randomly initialized systems settle into a steady state crystal lattice. We compute that all LCE's become negative in this region, and this is the only region that is not chaotic. 

Slightly lower values of $\mu_y$ will still produce a global stable crystal lattice of cells with high $x$/low $y$ values, but with chaotic waves washing through the surrounding low $x$/high $y$ cells, as shown in the middle row of figure \ref{fig:crystals}. This region produces positive Lyapunov exponents.    

As the $\mu_y$ value is further lowered, the chaotic waves get more pronounced and break up the crystalline pattern into transient islands of crystal structures, shown in the bottom row of figure \ref{fig:crystals}. These organized islands break apart as waves pass through and new islands are formed in the wake.

\subsection{Gingham Dynamics}\label{sec:gingham}

\begin{figure}[!htbp]
\begin{center}
  \includegraphics[width=.45\linewidth]{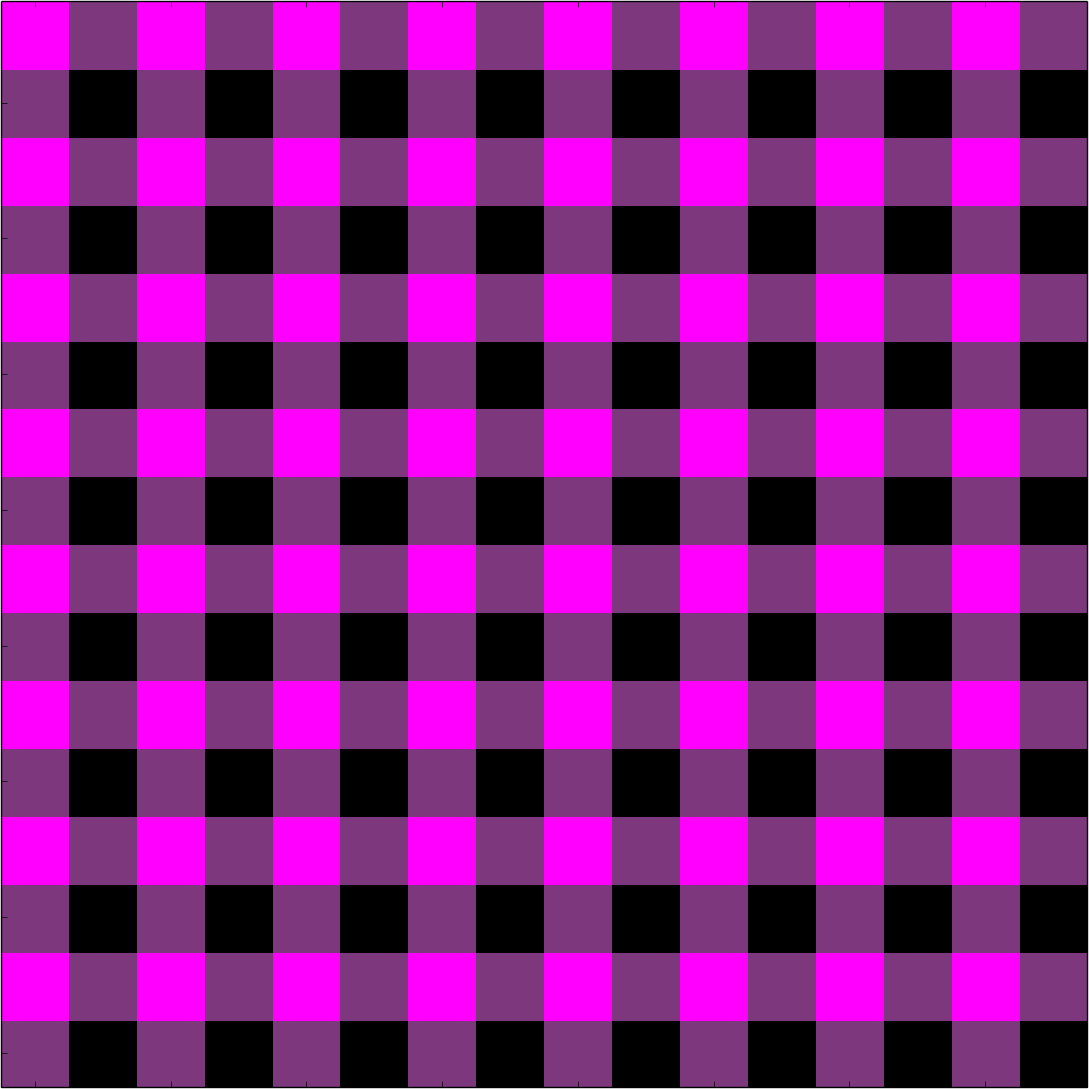}
  \hskip6pt
  \includegraphics[width=.45\linewidth]{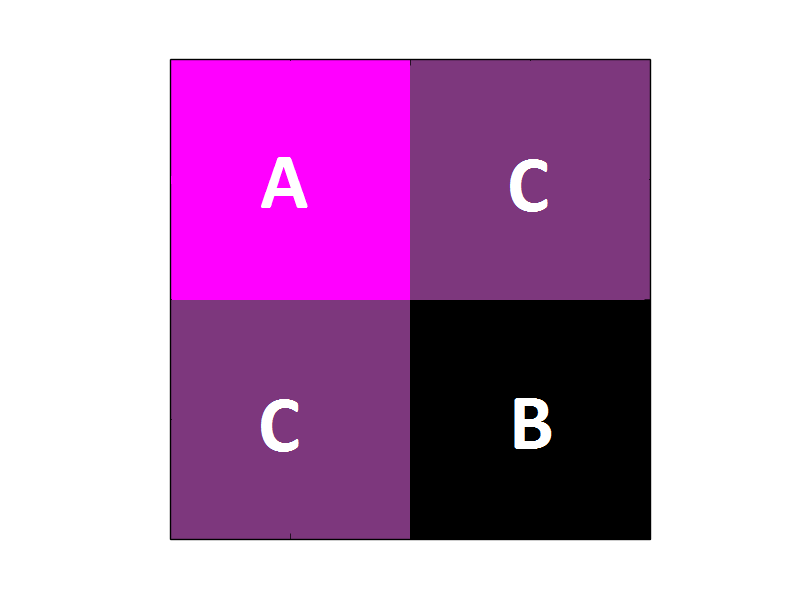}
\caption{Uniform crystalline structure and fundamental domain.}
\label{fig:perfect}
\end{center}
\end{figure}

A perfect crystalline pattern is shown in figure \ref{fig:perfect}, which is colored according to the host/prey population ($x$) and shows a clear gingham pattern. This is seen to be tessellated by a smaller $2\times 2$ grid modeled with three pairs of variables. 

Applying the Nicholson-Bailey dynamics to this configuration results in a six-dimensional system.

$$
\begin{array}{rcl}
x_A^{t+1}  &=& (1-\mu_x)f(x_A^t,y_A^t) + \textstyle{1\over 2}\mu_x f(x_B^t,y_B^t) + \textstyle{1\over 2}\mu_x f(x_C^t,y_C^t) \\
y_A^{t+1}  &=& (1-\mu_y)g(x_A^t,y_A^t) + \textstyle{1\over 2}\mu_y g(x_B^t,y_B^t) + \textstyle{1\over 2}\mu_y g(x_C^t,y_C^t) \\
x_B^{t+1}  &=& (1-\mu_x)f(x_B^t,y_B^t) + \textstyle{1\over 2}\mu_x f(x_A^t,y_A^t) + \textstyle{1\over 2}\mu_x f(x_C^t,y_C^t) \\
y_B^{t+1}  &=& (1-\mu_y)g(x_B^t,y_B^t) + \textstyle{1\over 2}\mu_y g(x_A^t,y_A^t) + \textstyle{1\over 2}\mu_y g(x_C^t,y_C^t) \\
x_C^{t+1}  &=& (1-\mu_x)f(x_C^t,y_C^t) + \textstyle{1\over 4}\mu_x f(x_A^t,y_A^t) + \textstyle{1\over 4}\mu_x f(x_B^t,y_B^t) + \textstyle{1\over 2}\mu_x f(x_C^t,y_C^t) \\
y_C^{t+1}  &=& (1-\mu_y)g(x_C^t,y_C^t) + \textstyle{1\over 4}\mu_y g(x_A^t,y_A^t) + \textstyle{1\over 4}\mu_y g(x_B^t,y_B^t) + \textstyle{1\over 2}\mu_y g(x_C^t,y_C^t) \\
\end{array}
$$

\begin{figure}[!htbp]
\begin{center}
  \includegraphics[width=.45\linewidth]{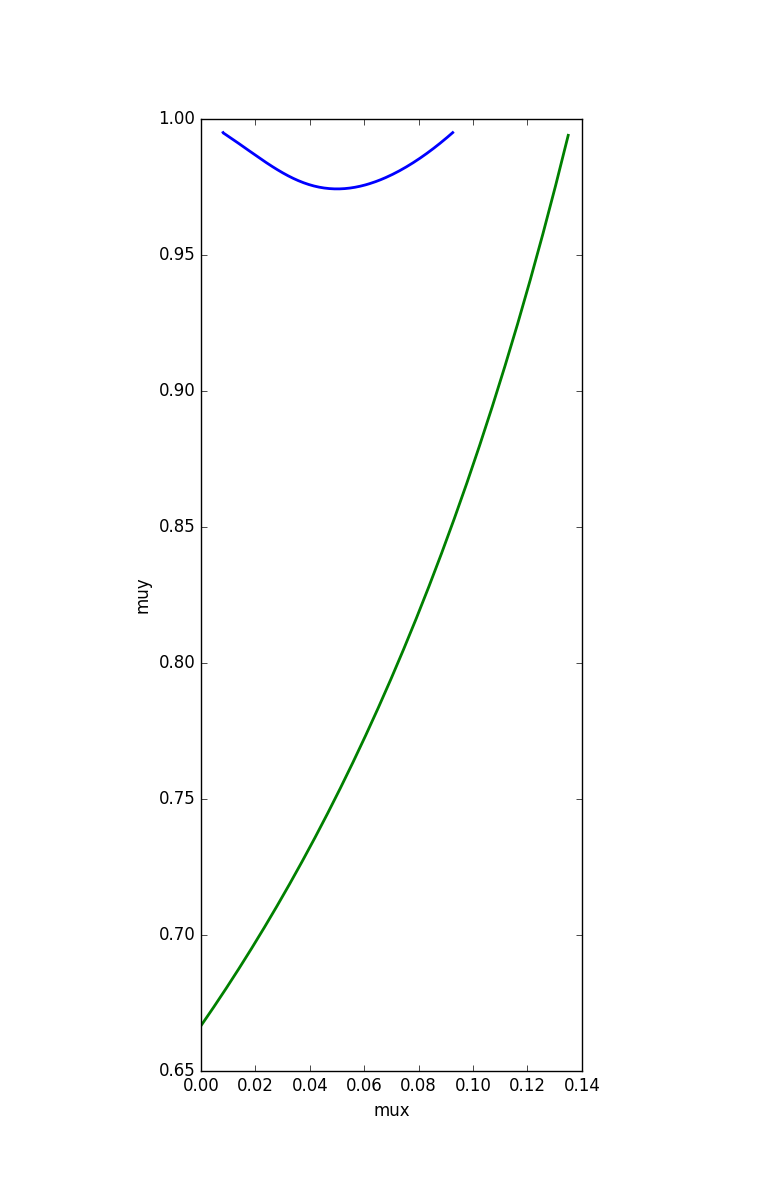}
  \includegraphics[width=.45\linewidth]{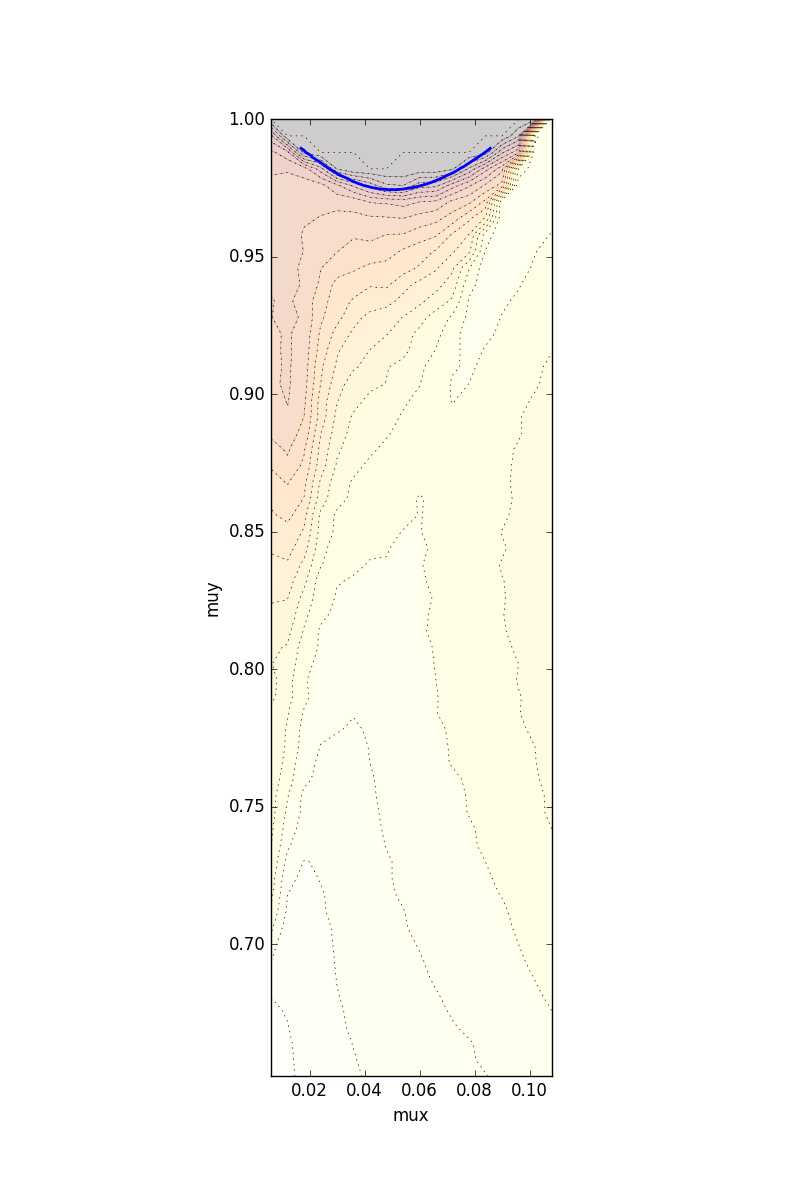}
\caption{Bifurcations of symmetric fixed point in reduced system.}
\label{fig:bifs}
\end{center}
\end{figure}

This system always has a fixed point at the the Nicholson-Bailey fixed point (\ref{eqn:fp}) with $x_A=x_B=x_C=x^*$, $y_A=y_B=y_C=y^*$. This fixed point bifurcates along the green curve in figure \ref{fig:bifs} to a symmetric pair of unstable fixed points. These fixed points show the typical crystal pattern where either the $A$ or $B$ cell has high $x$ / low $y$ values, and the remaining cells have low $x$ / high $y$ values. For example, at $\mu_x=0.05$, $\mu_y=0.99$ we find the stable fixed point at approximately 
$$
\begin{array}{lclcl}
x_A= 26.17 &\quad& x_B=0.69 &\quad& x_C= 0.37 \\
y_A=\ 0.64   &\quad& y_B=6.32 &\quad& y_C=3.42
\end{array}  
$$

The symmetric pair of fixed points become stable above the blue curve in figure \ref{fig:bifs}, synonymous with the region for negative Lyapunov exponents in the larger Nicholson-Bailey dynamic, providing excellent insight into the stable crystalline behavior.     

\section{Transitions and Spatial Chimeras}\label{sec:trans}

Lyapunov sampling relies on the assumption that (1) the system is spatially homogeneous and (2) convergence occurs on a reasonable computational timescale. Homogeneity can be tested by comparing distributions of LCE's across sampling patches and convergence can be monitored -- see Appendix II. 

There are two situations that violate one or both of these assumptions: behavior near bifurcations, and spatial chimeras.   

\subsection{Bifurcations}

As parameter values approach a bifurcation curve, some post-bifurcation behaviors may show up in spatially limited areas and persist for long periods in the system. When we use Lyapunov sampling, the samples in these regions will produce a Lyapunov spectrum more characteristic of the post-bifurcation region, where other samples will reflect the pre-bifurcation spectrum.  

\begin{figure}[!htbp]
\begin{center}
  \includegraphics[width=.4\linewidth, height=.6\linewidth]{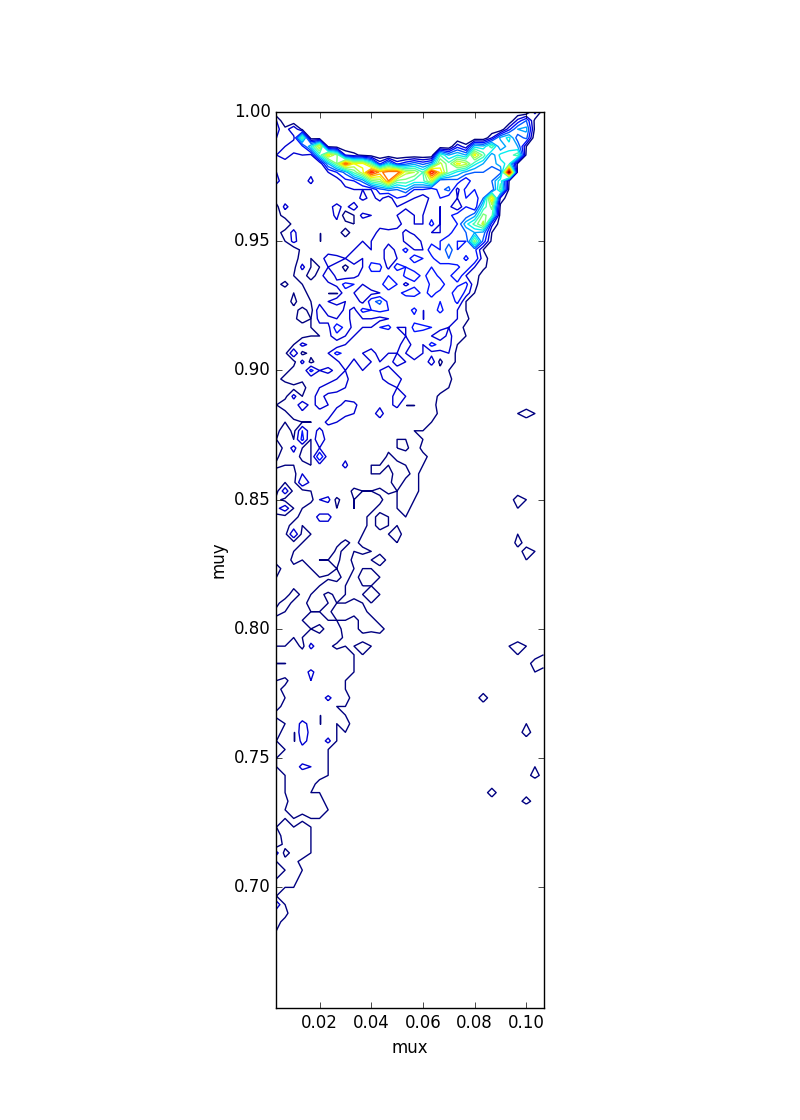}
  \includegraphics[width=.3\linewidth, height=.6\linewidth]{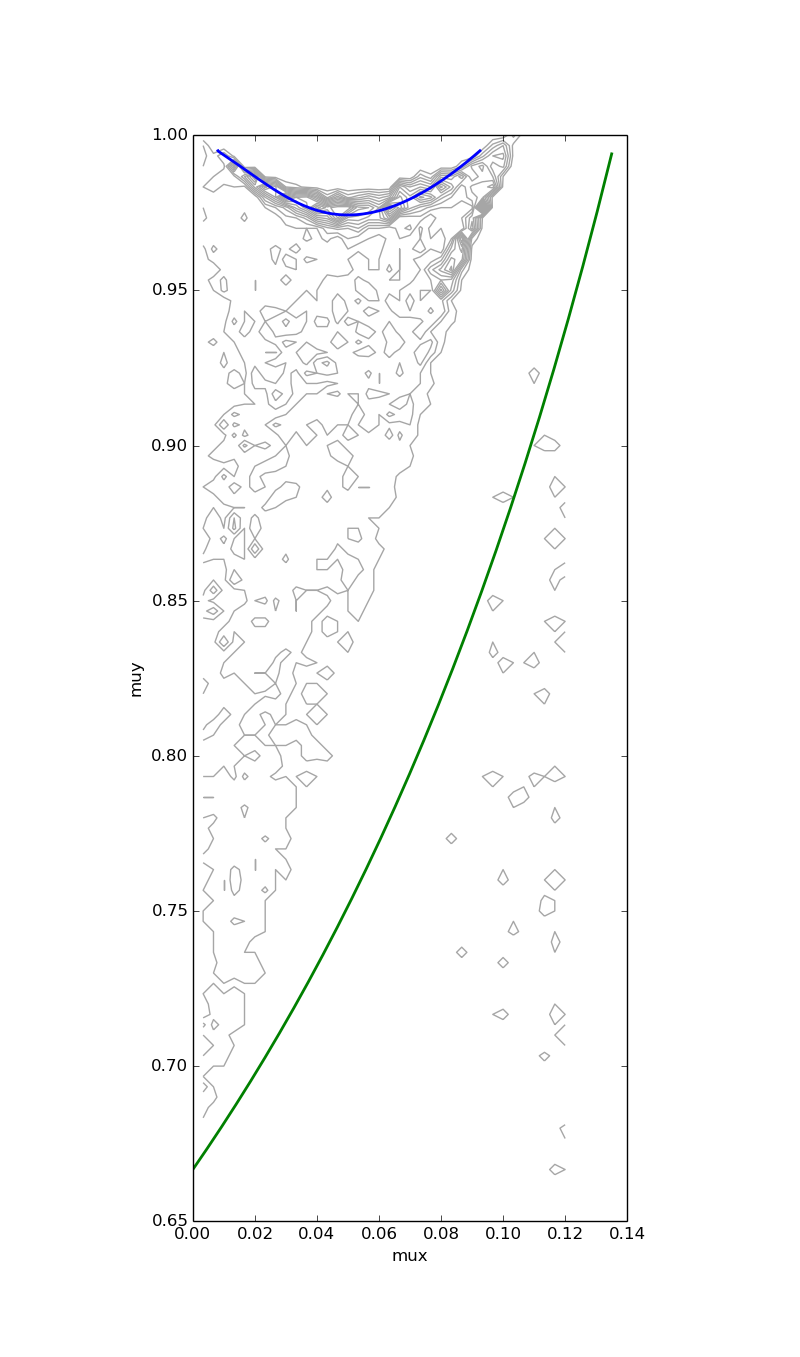}
\caption{Contour of difference in MLE's for $\lambda=2$, small $\mu_x$, and large $\mu_y$.}
\label{fig:NWsamp}
\end{center}
\end{figure}

For example, in section \ref{sec:crystal} we identified a semi-circular region for small $\mu_x$ and large $\mu_y$ with no positive Lyapunov exponents. Near the boundary of this region Lyapunov sampling produces a wide discrepancy in the maximum Lyapunov exponents (MLE) between subgrids. Some regions of the randomly seeded plane quickly settle into a stationary crystal pattern and a sample subgrid in this area will produce negative exponents. But because the system is near a bifurcation other areas remain active for long periods of time, producing positive exponents in the calculation. Thus there is a high discrepancy in the MLE between the samples. Figure \ref{fig:NWsamp} is contour plot of the difference between the largest MLE and smallest MLE of the sampled regions, with very high differences near the boundary of the semi-circular region. We used a $768 \times 512$ grid with an array of sample subgrids of size $16\times16$, for which the LCE spectrum was computed over $6000$ iterates. 

This figure also indicates a remarkably linear boundary creating a triangular region in the northwest corner. The transition between a global stable crystal lattice and transient crystalline islands is gradual, no clear boundary is evident in the data. However, the boundary between transient crystalline islands and no crystalline patterns at all seems to occur along this linear boundary. The crystalline islands can be long-lived and will depress the computed exponents in a sample that is dominated by an island, thus leading to the discrepancy in MLE's shown in figure \ref{fig:NWsamp}. The right-hand graph in figure \ref{fig:NWsamp} shows the MLE contours in relation to the bifurcation curves discussed in section \ref{sec:gingham}.

\subsection{Toroidal Chimeras} 

Near the $\mu_x=1$, $\mu_y=0$ corner, and extending along the $\mu_y=0$ boundary, we find coexisting regions of chaos and organized waves. 

\begin{figure}[!htbp]
	\begin{center}
		\includegraphics[width=.8\linewidth]{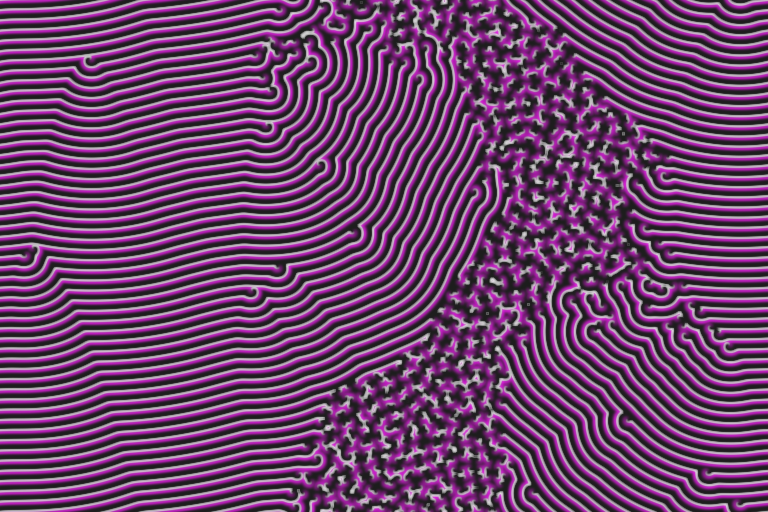}
		\caption{Emergent chimera $\lambda=2$ on a $768 \times 512$ grid, $\mu_x=0.95$, $\mu_y=0.05$.  }
		\label{fig:chim}
	\end{center}
\end{figure}

For example, $\mu_x=0.95$ and $\mu_y=0.05$ produces a chimeric structure of highly organized areas of waves cut through by a swath of more chaotic dynamics, as shown in figure \ref{fig:chim}. The chaotic swath favors the shorter distance around the torus, and slowly rotates around the longer direction over a timescale of millions of iterates.

\begin{figure}[!htbp]
\begin{center}
  \includegraphics[width=.6\linewidth]{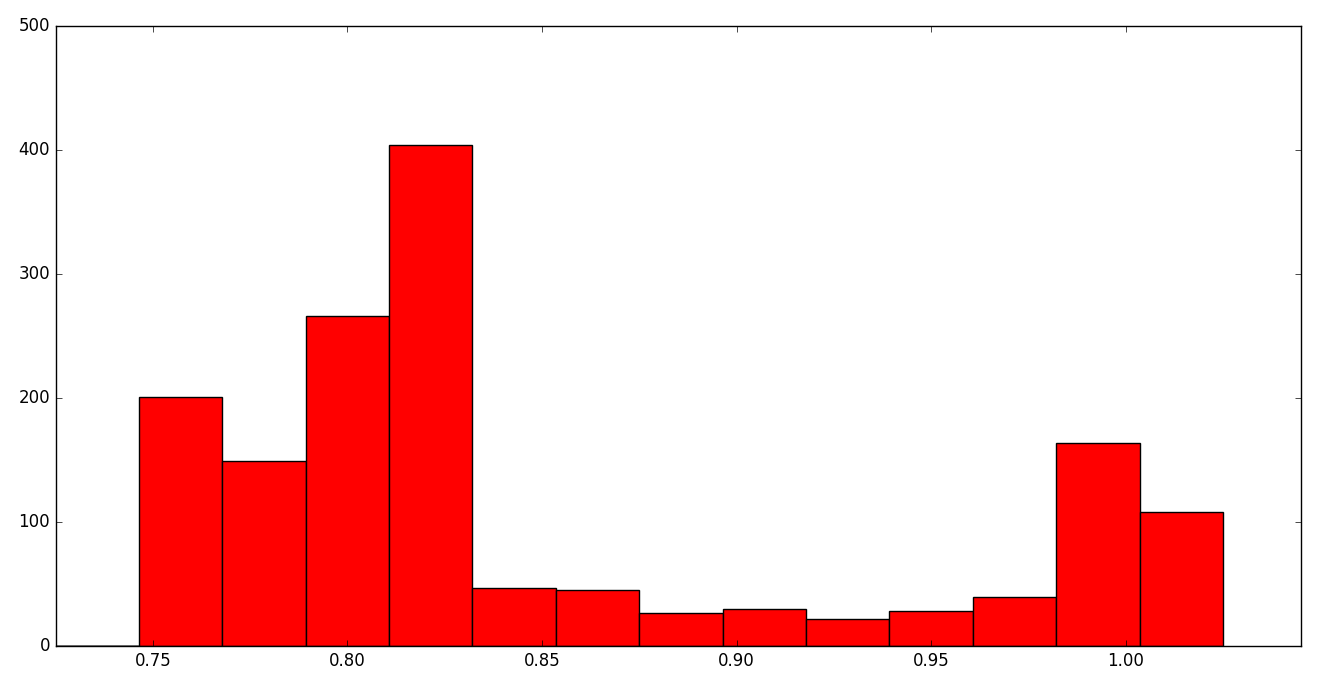}
\caption{Distribution of maximal Lyapunov exponents from toroidal chimera.}
\label{fig:bimod}
\end{center}
\end{figure}

Figure \ref{fig:chim} shows the dynamics on a $768 \times 512$ grid after 20 million iterates from a random initialization. The grid is tessellated with a $48\times32$ array of sample subgrids of size $16\times16$, for which the LCE spectrum was computed over $6000$ iterates. The distribution of maximum exponents from the samples demonstrates clear bimodality, shown in figure \ref{fig:bimod}, due to some samples being in the chaotic region and others in the organized region.

\section{Conclusions and Thoughts}\label{sec:Fin}

Lyapunov sampling is effective in approximating the Lyapunov spectrum in spatially homogeneous systems. Sampling can also be used to indicate the presence of emergent structures and delineate bifurcations. 

The region where crystalline islands emerge seems to have a remarkably linear boundary, suggesting some underlying phenomena that allows the persistence of these structures that would be worthwhile to try and understand. 

The computing platform for this work is unremarkable; using an NVIDIA GTX 1080 GPU on a host system with an Intel i7-6700K CPU @ 4.00GHz. Software tools consisted of Nvidia's Cuda C compiler and EM Photonics  numerics package CULA Tools.

\section*{Appendix I: Computing Lyapunov Exponents }\label{sec:AppendixI}

Understanding the theory of Lyapunov exponents requires consideration of hyperbolicity for generic points of an invariant measure, and a full treatment can be found in \cite{Pesin}. Without getting too deep into the abstract details, the basic ideas behind the numerical computation of Lyapunov exponents are as follows.  

Lyapunov exponents characterize the exponential rates at which trajectories become separated in an attractor. For a typical $x$ in the basin of the attractor, the Maximal Lyapunov Exponent (MLE) is the long-run exponential rate at which trajectories of generic points near $x$ in the same basin will become separated from the trajectory of $x$. For a transformation $T:\R^N\to\R^N$ the Jacobian of $T^n$ evaluated at a point $x$ is given as $\J(T^n)|_x$.  The MLE is given by  
$$\lambda= \max_{\|v\|=1} \left(\limsup_{n\to\infty} {\textstyle{1\over n}} \log \| \J(T^n)|_x\,v\|\right)$$ 
where $v\in \R^N$. The Jacobian expands by the product rule to a co-cycle: 
$$\J(T^n)|_x=\J(T)|_{T^{n-1}(x)} \; \J(T)|_{T^{n-2}(x)}\; \cdots \J(T)|_x.$$ 

Taking $J_k=\J(T)|_{T^{k-1}(x)}$ and defining $\chi(v) = \limsup_{n\to\infty} {1\over n} \log \|\prod_{k=1}^n J_k\; v\|$, the MLE is then 
$$\lambda= \max_{\|v\|=1} \chi(v).$$ 

Numerical computation of long products of matrices is notoriously unstable, and we use $QR$-factorization (decomposing a matrix into a product of an orthogonal matrix, $Q$, and upper triangular matrix, $R$) to preserve structure.

For $\prod J_k$, let $B_1=J_1$, and recursively define $Q_iR_i$ as the $QR$-factorization of $B_i$ with  $B_{i+1}=J_{i+1}Q_i$:

\renewcommand{\arraystretch}{1.0}
$$ \begin{array}{rcl}
\ldots J_3J_2J_1 &\quad&\\
& & J_1\to Q_1R_1 \\
\ldots J_3J_2Q_1R_1 & & \\
& & J_2Q_1 \to B_2 \\
\ldots J_3B_2R_1 & & \\
& & B_2 \to Q_2R_2 \\
\ldots J_3 Q_2R_2 R_1 \\
& &  J_3 Q_2 \to B_3\\
\hfil\vdots\hfil & & \hfil\vdots\hfil
\end{array}
$$

At the $n$th stage we have $J_n\cdots J_1 = Q_nR_n\cdots R_1$ with $Q_n$ unitary and $S_n= R_n\cdots R_1$ upper triangular. For the MLE, it is readily shown that  
$$\lambda= \max_{\|v\|=1} \chi(v)= \max_i \left(\limsup_{n\to\infty} {\textstyle{1\over n}}\log\|S_n[i,i]\|\right)$$
using basis decomposition and the properties of $\chi(\cdot)$. Recursively removing the maximizing direction and repeating produces the full Lyapunov spectrum.

In fact, for well-behaved uniformly hyperbolic systems the exponential rate of growth of the individual diagonal elements of $S_n$ converge to the Lyapunov exponents. That is, the spectrum of characteristic Lyapunov exponents is given by $$\lim_{n\to\infty} {\textstyle{1\over n}}\log\|S_n[i,i]\|$$ for $i=1,\ldots,N$. This convergence is usually assumed when computing Lyapunov exponents, and the assumption can be supported by examining the numerical convergence.  

\section*{Appendix II: Subgrid Sampling }\label{sec:AppendixII}

To demonstrate the viability of subgrid sampling for estimating Lyapunov exponents in the Nicholson-Bailey mode, we take take $\lambda=2$, $\mu_x=0.6$ and $\mu_y=0.8$ and compute the complete LCE spectrum for a $96 \times 64$ system as well as the spectra for three $32\times 32$ subgrids. 

\begin{figure}[!htbp]
\begin{center}
  \includegraphics[width=\linewidth]{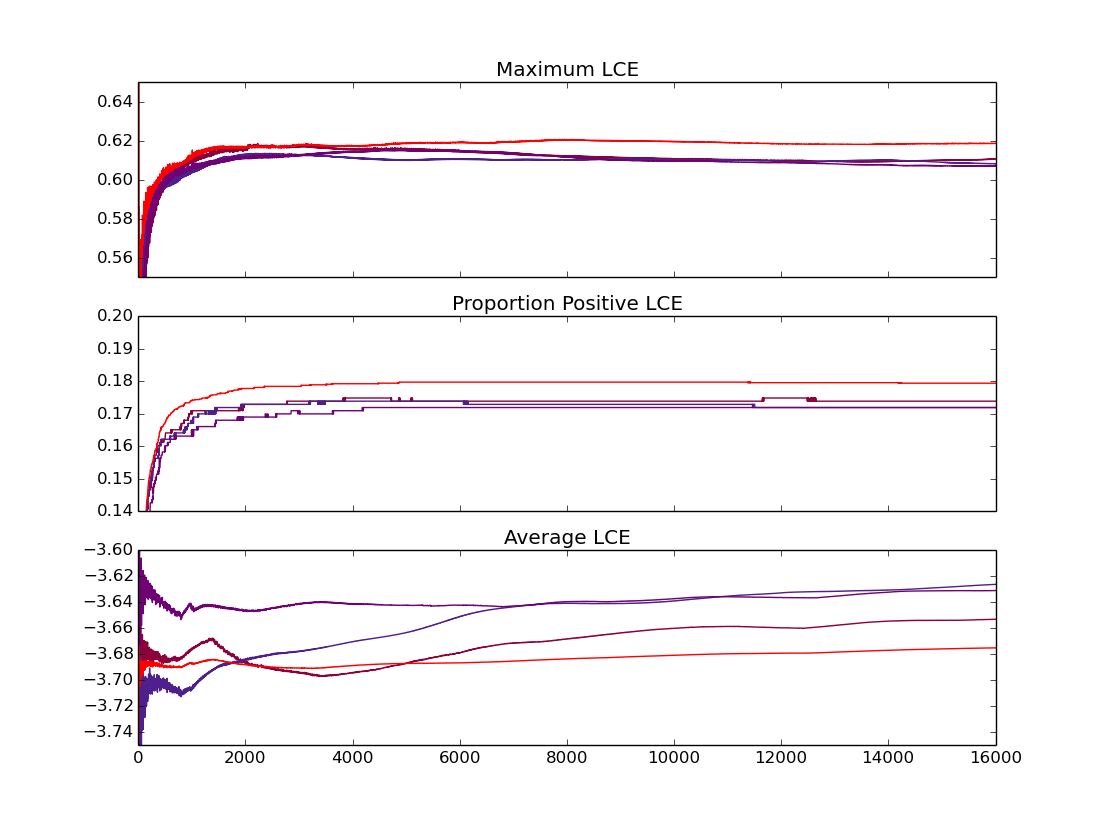}
\caption{Convergence of maximum LCE, proportion of positive LCE's, and average LCE for a larger grid (red) and three subgrids (purples).}
\label{fig:converge}
\end{center}
\end{figure}

\begin{figure}[!htbp]
\begin{center}
  \includegraphics[width=\linewidth]{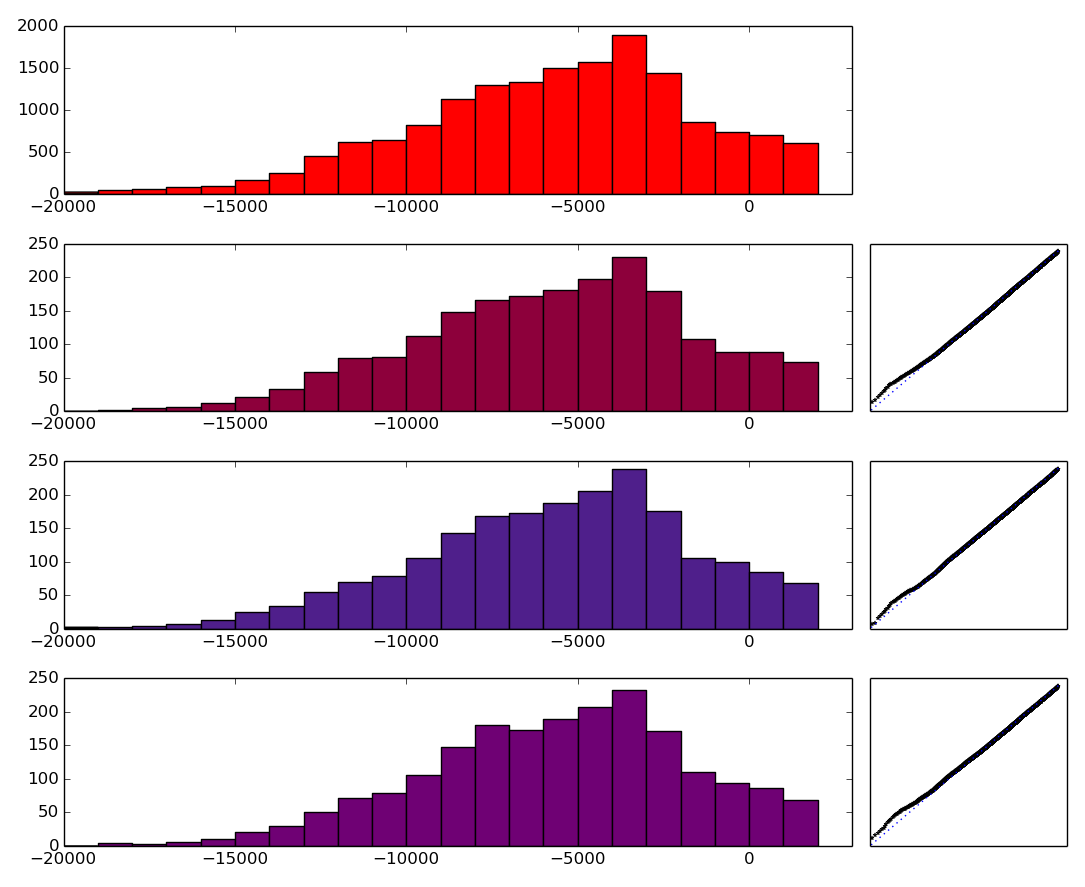}
\caption{Distribution of LCE's for larger grid (top) and three subgrids, with quantile plots.}
\label{fig:histograms}
\end{center}
\end{figure}

Convergence is shown in figure \ref{fig:converge}; the values for the maximum LCE and the proportion of positive LCE's converge fairly quickly and are reasonable approximations to the LCE's of the larger system. Of particular note is that the distribution of LCE's after 3000 iterates is remarkably similar for the subgrids and the larger grid (see figure \ref{fig:histograms}). These results are consistent with tests using larger systems (up to $96 \times 128$) and smaller subgrids ($16\times 16$). The subgrid sampling technique provides highly persuasive evidence for Rand's criteria for chaos in high dimensional homogeneous systems. 

\bibliographystyle{plain}

\bibliography{references}

\begin{thebibliography}{10}

\bibitem{SpatEco}
{\em Spatial Ecology: The Role of Space in Population Dynamics and
  Interspecific Interactions (MPB-30)}.
\newblock Princeton University Press, 1997.

\bibitem{Pesin}
L~Barreira and Y~Pesin.
\newblock {\em Nonuniform Hyperbolicity}.
\newblock Cambridge, 2007.

\bibitem{OttarJordi}
Ottar~N. Bjørnstad and Jordi Bascompte.
\newblock Synchrony and second-order spatial correlation in host–parasitoid
  systems.
\newblock {\em Journal of Animal Ecology}, 70(6):924--933, 2001.

\bibitem{CHM92}
HN~Comins, MP~Hassell, and RM~May.
\newblock {The spatial dynamics of host parasitoid systems}.
\newblock {\em {Journal of Animal Ecology}}, {61}({3}):{735--748}, {1992}.

\bibitem{CHM91}
MP~Hassell, HN~Comins, and RM~May.
\newblock {Spatial structure and chaos in insect population-dynamics}.
\newblock {\em {Nature}}, {353}({6341}):{255--258}, {SEP 19} {1991}.

\bibitem{CHM94}
MP~Hassell, HN~Comins, and RM~May.
\newblock {Species coexistence and self-organizing spatial dynamics}.
\newblock {\em {NATURE}}, {370}({6487}):{290--292}, {JUL 28} {1994}.

\bibitem{May95}
Robert~M. May.
\newblock Necessity and chance: deterministic chaos in ecology and evolution.
\newblock {\em Bull. Amer. Math. Soc. (N.S.)}, 32(3):291--308, 1995.

\bibitem{Nich33}
A.J. Nicholson.
\newblock {The balance of animal populations}.
\newblock {\em {Journal of Animal Ecology}}, {2}:{132--178}, {1933}.

\bibitem{Nich35}
A.J. Nicholson and V.A. Bailey.
\newblock {The balance of animal populations, Part I}.
\newblock {\em {Proceedings of the Zoological Society of London}},
  {3}:{551--598}, {1935}.

\bibitem{Rand94}
DA~Rand.
\newblock {Measuring and characterizing spatial patterns, dynamics and chaos in
  spatially extended dynamical-systems and ecologies}.
\newblock {\em {Philosophical Transactions of the Royal Society A-Mathematical
  Physical and Engineering Sciences}}, {348}({1688}):{497--514}, {SEP 15}
  {1994}.

\bibitem{SoleValls}
RV~Sole, J~Valls, and J~Bascompte.
\newblock {Spiral waves, chaos and multiple attractors in lattice models of
  interacting populations}.
\newblock {\em {Physics Letters A}}, {166}({2}):{123--128}, {JUN 15} {1992}.

\bibitem{Hub1}
Hubertus~F. von Bremen, Firdaus~E. Udwadia, and Wlodek Proskurowski.
\newblock An efficient {QR} based method for the computation of lyapunov
  exponents.
\newblock {\em Phys. D}, 101(1-2):1--16, 1997.

\bibitem{WakanoHauert}
Joe~Yuichiro Wakano and Christoph Hauert.
\newblock {Pattern formation and chaos in spatial ecological public goods
  games}.
\newblock {\em {Journal of Theoretical Biology}}, {268}({1}):{30--38}, {JAN 7}
  {2011}.

\end{thebibliography}

\end{document}